\documentclass[12pt]{article}
\usepackage{graphicx}
\usepackage{xcolor}
\usepackage{amssymb,amsfonts,amsthm,amsmath,calligra,tikz}
\usepackage[utf8]{inputenc}
\usepackage{accents}
\usepackage{slashed}
\usepackage{yfonts}
\usepackage{mathrsfs,pifont}
\theoremstyle{definition}
\usepackage{tikz}
\usepackage[all]{xy}

\def\be{\begin{eqnarray}}
\def\ee{\end{eqnarray}}

\def\matZ{{\mathbb{Z}}}
\def\matR{{\mathbb{R}}}

\def\matC{{\mathbb{C}}}

\newcommand{\bA}{\mathsf{A}}

\newcommand{\bT}{\mathsf{T}}
\newcommand{\bK}{\mathsf{K}}

\newcommand{\Lie}{\mathrm{Lie}}
\newcommand{\Pic}{\mathrm{Pic}}

\newcommand{\cL}{\mathscr{L}}
\newcommand{\cP}{\mathscr{P}}
\let\bs\boldsymbol

\def\zz{{\bs z}}

\def\aa{{\bs a}}

\theoremstyle{definition}

\newtheorem{Conjecture}{Conjecture}
\newtheorem{Theorem}{Theorem}

\newtheorem{Note}{Note}

\newcommand{\hilb}{X}

\newcommand{\fC}{\mathfrak{C}} 

\def\tb {{\cal{V}}}  

\makeatletter
\newcommand{\somespecialrotate}[3][]{%
\begingroup
\sbox\@tempboxa{#3}%
\@tempdima=.5\wd\@tempboxa
\sbox\@tempboxa{\rotatebox[#1]{#2}{\usebox\@tempboxa}}%
\advance\@tempdima by -.5\wd\@tempboxa
\mbox{\hskip\@tempdima\usebox\@tempboxa}%
\endgroup}
\linethickness{1 pt}
\newsavebox{\ybb}
\savebox{\ybb}(1,1){
	\put(0, 0){ \color{blue} \line(0, 1){ 36}} \put(2, 0){ \color{blue} \line(0, 1){ 36}} \put(4, 0){ \color{blue} \line(0, 1){ 30}} \put(6, 0){ \color{blue} \line(0, 1){ 28}} \put(8, 0){ \color{blue} \line(0, 1){ 26}} \put(10, 0){ \color{blue} \line(0, 1){ 24}} \put(12, 0){ \color{blue} \line(0, 1){ 20}} \put(14, 0){ \color{blue} \line(0, 1){ 16}} \put(16, 0){ \color{blue} \line(0, 1){ 16}} \put(18, 0){ \color{blue} \line(0, 1){ 14}} \put(20, 0){ \color{blue} \line(0, 1){ 12}} \put(22, 0){ \color{blue} \line(0, 1){ 10}} \put(24, 0){ \color{blue} \line(0, 1){ 10}} \put(26, 0){ \color{blue} \line(0, 1){ 8}} \put(28, 0){ \color{blue} \line(0, 1){ 6}} \put(30, 0){ \color{blue} \line(0, 1){ 4}} \put(32, 0){ \color{blue} \line(0, 1){ 2}} \put(34, 0){ \color{blue} \line(0, 1){ 2}} \put(36, 0){ \color{blue} \line(0, 1){ 2}} \put(0, 0){ \color{blue} \line(1, 0){ 36}} \put(0, 2){ \color{blue} \line(1, 0){ 36}} \put(0, 4){ \color{blue} \line(1, 0){ 30}} \put(0, 6){ \color{blue} \line(1, 0){ 28}} \put(0, 8){ \color{blue} \line(1, 0){ 26}} \put(0, 10){ \color{blue} \line(1, 0){ 24}} \put(0, 12){ \color{blue} \line(1, 0){ 20}} \put(0, 14){ \color{blue} \line(1, 0){ 18}} \put(0, 16){ \color{blue} \line(1, 0){ 16}} \put(0, 18){ \color{blue} \line(1, 0){ 12}} \put(0, 20){ \color{blue} \line(1, 0){ 12}} \put(0, 22){ \color{blue} \line(1, 0){ 10}} \put(0, 24){ \color{blue} \line(1, 0){ 10}} \put(0, 26){ \color{blue} \line(1, 0){ 8}} \put(0, 28){ \color{blue} \line(1, 0){ 6}} \put(0, 30){ \color{blue} \line(1, 0){ 4}} \put(0, 32){ \color{blue} \line(1, 0){ 2}} \put(0, 34){ \color{blue} \line(1, 0){ 2}} \put(0, 36){ \color{blue} \line(1, 0){ 2}}
	\put(1, 1){ \color{red} \line(0, 1){ 2}} \put(1, 3){ \color{red} \line(0, 1){ 2}} \put(1, 5){ \color{red} \line(0, 1){ 2}} \put(1, 7){ \color{red} \line(0, 1){ 2}} \put(1, 7){ \color{red} \line(1, 0){ 2}} \put(1, 9){ \color{red} \line(0, 1){ 2}} \put(1, 9){ \color{red} \line(1, 0){ 2}} \put(1, 11){ \color{red} \line(0, 1){ 2}} \put(1, 11){ \color{red} \line(1, 0){ 2}} \put(1, 13){ \color{red} \line(0, 1){ 2}} \put(1, 13){ \color{red} \line(1, 0){ 2}} \put(1, 15){ \color{red} \line(0, 1){ 2}} \put(1, 15){ \color{red} \line(1, 0){ 2}} \put(1, 17){ \color{red} \line(0, 1){ 2}} \put(1, 19){ \color{red} \line(0, 1){ 2}} \put(1, 19){ \color{red} \line(1, 0){ 2}} \put(1, 21){ \color{red} \line(0, 1){ 2}} \put(1, 21){ \color{red} \line(1, 0){ 2}} \put(1, 23){ \color{red} \line(0, 1){ 2}} \put(1, 23){ \color{red} \line(1, 0){ 2}} \put(1, 25){ \color{red} \line(0, 1){ 2}} \put(1, 25){ \color{red} \line(1, 0){ 2}} \put(1, 27){ \color{red} \line(0, 1){ 2}} \put(1, 27){ \color{red} \line(1, 0){ 2}} \put(1, 29){ \color{red} \line(0, 1){ 2}} \put(1, 29){ \color{red} \line(1, 0){ 2}} \put(1, 31){ \color{red} \line(0, 1){ 2}} \put(1, 33){ \color{red} \line(0, 1){ 2}} \put(3, 1){ \color{red} \line(0, 1){ 2}} \put(3, 1){ \color{red} \line(1, 0){ 2}} \put(3, 3){ \color{red} \line(0, 1){ 2}} \put(3, 3){ \color{red} \line(1, 0){ 2}} \put(3, 5){ \color{red} \line(0, 1){ 2}} \put(3, 5){ \color{red} \line(1, 0){ 2}} \put(3, 15){ \color{red} \line(1, 0){ 2}} \put(3, 17){ \color{red} \line(0, 1){ 2}} \put(3, 21){ \color{red} \line(1, 0){ 2}} \put(3, 23){ \color{red} \line(1, 0){ 2}} \put(3, 25){ \color{red} \line(1, 0){ 2}} \put(3, 27){ \color{red} \line(1, 0){ 2}} \put(5, 1){ \color{red} \line(1, 0){ 2}} \put(5, 3){ \color{red} \line(1, 0){ 2}} \put(5, 7){ \color{red} \line(0, 1){ 2}} \put(5, 7){ \color{red} \line(1, 0){ 2}} \put(5, 9){ \color{red} \line(0, 1){ 2}} \put(5, 11){ \color{red} \line(0, 1){ 2}} \put(5, 11){ \color{red} \line(1, 0){ 2}} \put(5, 13){ \color{red} \line(0, 1){ 2}} \put(5, 13){ \color{red} \line(1, 0){ 2}} \put(5, 17){ \color{red} \line(0, 1){ 2}} \put(5, 19){ \color{red} \line(0, 1){ 2}} \put(5, 19){ \color{red} \line(1, 0){ 2}} \put(5, 21){ \color{red} \line(1, 0){ 2}} \put(5, 25){ \color{red} \line(1, 0){ 2}} \put(7, 5){ \color{red} \line(0, 1){ 2}} \put(7, 7){ \color{red} \line(1, 0){ 2}} \put(7, 9){ \color{red} \line(0, 1){ 2}} \put(7, 11){ \color{red} \line(1, 0){ 2}} \put(7, 15){ \color{red} \line(0, 1){ 2}} \put(7, 17){ \color{red} \line(0, 1){ 2}} \put(7, 17){ \color{red} \line(1, 0){ 2}} \put(7, 23){ \color{red} \line(0, 1){ 2}} \put(7, 23){ \color{red} \line(1, 0){ 2}} \put(9, 1){ \color{red} \line(0, 1){ 2}} \put(9, 3){ \color{red} \line(0, 1){ 2}} \put(9, 3){ \color{red} \line(1, 0){ 2}} \put(9, 5){ \color{red} \line(0, 1){ 2}} \put(9, 9){ \color{red} \line(0, 1){ 2}} \put(9, 9){ \color{red} \line(1, 0){ 2}} \put(9, 11){ \color{red} \line(1, 0){ 2}} \put(9, 13){ \color{red} \line(0, 1){ 2}} \put(9, 15){ \color{red} \line(0, 1){ 2}} \put(9, 19){ \color{red} \line(0, 1){ 2}} \put(9, 19){ \color{red} \line(1, 0){ 2}} \put(9, 21){ \color{red} \line(0, 1){ 2}} \put(11, 1){ \color{red} \line(0, 1){ 2}} \put(11, 3){ \color{red} \line(1, 0){ 2}} \put(11, 5){ \color{red} \line(0, 1){ 2}} \put(11, 5){ \color{red} \line(1, 0){ 2}} \put(11, 7){ \color{red} \line(0, 1){ 2}} \put(11, 7){ \color{red} \line(1, 0){ 2}} \put(11, 9){ \color{red} \line(1, 0){ 2}} \put(11, 13){ \color{red} \line(0, 1){ 2}} \put(11, 13){ \color{red} \line(1, 0){ 2}} \put(11, 15){ \color{red} \line(0, 1){ 2}} \put(11, 15){ \color{red} \line(1, 0){ 2}} \put(11, 17){ \color{red} \line(0, 1){ 2}} \put(13, 1){ \color{red} \line(0, 1){ 2}} \put(13, 1){ \color{red} \line(1, 0){ 2}} \put(13, 5){ \color{red} \line(1, 0){ 2}} \put(13, 7){ \color{red} \line(1, 0){ 2}} \put(13, 9){ \color{red} \line(1, 0){ 2}} \put(13, 11){ \color{red} \line(0, 1){ 2}} \put(13, 13){ \color{red} \line(1, 0){ 2}} \put(13, 15){ \color{red} \line(1, 0){ 2}} \put(15, 1){ \color{red} \line(1, 0){ 2}} \put(15, 3){ \color{red} \line(0, 1){ 2}} \put(15, 5){ \color{red} \line(1, 0){ 2}} \put(15, 7){ \color{red} \line(1, 0){ 2}} \put(15, 11){ \color{red} \line(0, 1){ 2}} \put(15, 11){ \color{red} \line(1, 0){ 2}} \put(15, 13){ \color{red} \line(1, 0){ 2}} \put(17, 3){ \color{red} \line(0, 1){ 2}} \put(17, 7){ \color{red} \line(1, 0){ 2}} \put(17, 9){ \color{red} \line(0, 1){ 2}} \put(17, 11){ \color{red} \line(1, 0){ 2}} \put(19, 1){ \color{red} \line(0, 1){ 2}} \put(19, 1){ \color{red} \line(1, 0){ 2}} \put(19, 3){ \color{red} \line(0, 1){ 2}} \put(19, 3){ \color{red} \line(1, 0){ 2}} \put(19, 5){ \color{red} \line(0, 1){ 2}} \put(19, 5){ \color{red} \line(1, 0){ 2}} \put(19, 7){ \color{red} \line(1, 0){ 2}} \put(19, 9){ \color{red} \line(0, 1){ 2}} \put(19, 9){ \color{red} \line(1, 0){ 2}} \put(21, 3){ \color{red} \line(1, 0){ 2}} \put(21, 5){ \color{red} \line(1, 0){ 2}} \put(21, 7){ \color{red} \line(1, 0){ 2}} \put(21, 9){ \color{red} \line(1, 0){ 2}} \put(23, 1){ \color{red} \line(0, 1){ 2}} \put(23, 5){ \color{red} \line(1, 0){ 2}} \put(23, 7){ \color{red} \line(1, 0){ 2}} \put(25, 1){ \color{red} \line(0, 1){ 2}} \put(25, 1){ \color{red} \line(1, 0){ 2}} \put(25, 3){ \color{red} \line(0, 1){ 2}} \put(25, 5){ \color{red} \line(1, 0){ 2}} \put(27, 3){ \color{red} \line(0, 1){ 2}} \put(27, 3){ \color{red} \line(1, 0){ 2}} \put(29, 1){ \color{red} \line(0, 1){ 2}} \put(29, 1){ \color{red} \line(1, 0){ 2}} \put(31, 1){ \color{red} \line(1, 0){ 2}} \put(33, 1){ \color{red} \line(1, 0){ 2}}}

\newsavebox{\ydb}
\savebox{\ydb}(1,1){ \put(4.5,2){\color{green}\rule{16pt}{16pt}}
	
	\put(0, 0){ \color{blue} \line(0, 1){ 10}} \put(2, 0){ \color{blue} \line(0, 1){ 10}} \put(4, 0){ \color{blue} \line(0, 1){ 6}} \put(6, 0){ \color{blue} \line(0,1){4}} \put(8, 0){ \color{blue} \line(0, 1){ 4}} \put(10, 0){ \color{blue} \line(0, 1){ 2}} \put(0, 0){ \color{blue} \line(1, 0){10}} \put(0, 2){ \color{blue} \line(1, 0){ 10}} \put(0, 4){ \color{blue} \line(1, 0){ 8}} \put(0, 6){ \color{blue} \line(1, 0){ 4}} \put(0, 8){ \color{blue} \line(1, 0){ 2}} \put(0, 10){ \color{blue} \line(1, 0){ 2}}
}

\newsavebox{\exoneb}
\savebox{\exoneb}(1,1){
		\put(0, 0){ \color{blue} \line(0, 1){ 4} } \put(2, 0){ \color{blue} \line(0, 1){ 4} } \put(4, 0){ \color{blue} \line(0, 1){ 4} } \put(0, 0){ \color{blue} \line(1, 0){ 4} } \put(0, 2){ \color{blue} \line(1, 0){ 4} } \put(0, 4){ \color{blue} \line(1, 0){ 4} }
	\put(1, 1){ \color{red} \line(0, 1){ 2}}
	\put(1, 1){ \color{red} \line(1, 0){ 2}}
	\put(1, 3){ \color{red} \line(1, 0){ 2}} }

\newsavebox{\extwob}
\savebox{\extwob}(1,1){
		\put(0, 0){ \color{blue} \line(0, 1){ 4} } \put(2, 0){ \color{blue} \line(0, 1){ 4} } \put(4, 0){ \color{blue} \line(0, 1){ 4} } \put(0, 0){ \color{blue} \line(1, 0){ 4} } \put(0, 2){ \color{blue} \line(1, 0){ 4} } \put(0, 4){ \color{blue} \line(1, 0){ 4} }
	\put(1, 1){ \color{red} \line(0, 1){ 2}}
	\put(3, 1){ \color{red} \line(0, 1){ 2}}
	\put(1, 3){ \color{red} \line(1, 0){ 2}}}

\newsavebox{\rgb}

\savebox{\rgb}(1,1)
{   \put(4.5,2){\color{blue}\rule{16pt}{16pt}} 
	\put(6.5,4){\color{green}\rule{16pt}{16pt}}
	\put(8.5,6){\color{green}\rule{16pt}{16pt}}
	\put(10.5,6){\color{red}\rule{16pt}{16pt}}
	\put(8.5,4){\color{red}\rule{16pt}{16pt}}
	\put(6.5,2){\color{red}\rule{16pt}{16pt}}
	\put(2.5,0){\color{red}\rule{16pt}{16pt}}
	\put(0.5,0){\color{green}\rule{16pt}{16pt}}
	\put(2.5,2){\color{green}\rule{16pt}{16pt}}
	\put(0, 0) { \color{blue} \line(0, 1) { 12} } \put(2, 0) { \color{blue} \line(0, 1) { 12} } \put(4, 0) { \color{blue} \line(0, 1) { 10} } \put(6, 0) { \color{blue} \line(0, 1) { 10} } \put(8, 0) { \color{blue} \line(0, 1) { 10} } \put(10, 0) { \color{blue} \line(0, 1) { 8} } \put(12, 0) { \color{blue} \line(0, 1) { 8} } \put(14, 0) { \color{blue} \line(0, 1) { 4} } \put(0, 0) { \color{blue} \line(1, 0) { 14} } \put(0, 2) { \color{blue} \line(1, 0) { 14} } \put(0, 4) { \color{blue} \line(1, 0) { 14} } \put(0, 6) { \color{blue} \line(1, 0) { 12} } \put(0, 8) { \color{blue} \line(1, 0) { 12} } \put(0, 10) { \color{blue} \line(1, 0) { 8} } \put(0, 12) { \color{blue} \line(1, 0) { 2} }}

\newsavebox{\ppb}

\savebox{\ppb}(1,1)
{   \put(0.5,10){\color{gray}\rule{16pt}{16pt}} 
	\put(2.5,8){\color{gray}\rule{16pt}{16pt}}
	 \put(0.5,8){\color{gray}\rule{16pt}{16pt}}
	 \put(0.5,6){\color{gray}\rule{16pt}{16pt}} 
	\put(0, 0) { \color{blue} \line(0, 1) { 12} } \put(2, 0) { \color{blue} \line(0, 1) { 12} } \put(4, 0) { \color{blue} \line(0, 1) { 10} } \put(6, 0) { \color{blue} \line(0, 1) { 10} } \put(8, 0) { \color{blue} \line(0, 1) { 10} } \put(10, 0) { \color{blue} \line(0, 1) { 8} } \put(12, 0) { \color{blue} \line(0, 1) { 8} } \put(14, 0) { \color{blue} \line(0, 1) { 4} } \put(0, 0) { \color{blue} \line(1, 0) { 14} } \put(0, 2) { \color{blue} \line(1, 0) { 14} } \put(0, 4) { \color{blue} \line(1, 0) { 14} } \put(0, 6) { \color{blue} \line(1, 0) { 12} } \put(0, 8) { \color{blue} \line(1, 0) { 12} } \put(0, 10) { \color{blue} \line(1, 0) { 8} } \put(0, 12) { \color{blue} \line(1, 0) { 2} }}

\begin{document}
\title{Characters of tangent spaces at torus fixed points and $3d$-mirror symmetry}
\author{Hunter Dinkins and Andrey Smirnov}
\date{}
\maketitle
\thispagestyle{empty}
	
\begin{abstract}
Let $X$ be a Nakajima quiver variety and $X'$ its $3d$-mirror. We consider the action of the Picard torus $\bK=\Pic(X)\otimes \matC^{\times}$ on $X'$.
Assuming that $(X')^{\bK}$ is finite, we propose a formula for the $\bK$-character of the tangent spaces at the fixed points in terms of certain enumerative invariants of~$X$ known as vertex functions.  
\end{abstract}
	
\setcounter{tocdepth}{2}

\section{Overview}

\subsection{} 
In this paper, we study symplectic varieties which appear as Higgs and Coulomb branches of certain  three-dimensional gauge theories with $\mathcal{N}=4$ supersymmetry. These theories were considered, for example, in \cite{BDG,GW,HW,SI}. 
We assume that the physical theories under consideration  are of quiver type, in which case the Higgs branches are known as {\it Nakajima quiver varieties}, see \cite{Nak1,GinzburgLectures} and Section 2 in \cite{MO} for an introduction\footnote{We believe that our main conjecture holds in full generality. We restrict the exposition to the quiver varieties for the sake of simplicity of the exposition.}. In this subsection we recall the most basic facts about these varieties important in the constructions below. 

The Nakajima varieties are smooth quasi-projective symplectic varieties equipped with a natural  action of an algebraic torus~$\bT$. The torus $\bT$ acts on a Nakajima variety $X$ by scaling the symplectic form
$\omega\in H^{2}(X,\matC)$. We denote by $\hbar \in \textrm{char}(\bT)$ the character of the one-dimensional $\bT$-module $\matC \omega$.  We denote by $\bA=\textrm{ker}(\hbar)\subset \bT$ the subtorus preserving the symplectic form and by $\matC^{\times}_{\hbar}:=\bT/\bA$ the corresponding one-dimensional factor.

The Nakajima quiver varieties are examples of symplectic resolutions and thus their cohomology are even \cite{Kaledin}.

The Nakajima varieties are defined as quotients by a group 
$$G=\prod\limits_{i\in I} GL(\textsf{v}_i)
$$
where $I$ denotes the (finite) set of vertices of the corresponding quiver. This means that 
$X$ is naturally equipped with a set of rank $\textsf{v}_i$ tautological vector 
bundles $\tb_{i}$. 

\begin{Theorem}[\cite{kirv}] \label{tmeven}
{\it 	If $X$ is a Nakajima variety then $K^{alg}_{\bT}(X)=K^{top}_{\bT}(X)$	is generated by Schur functors of tautological bundles $\tb_i$, $i\in I$. }
\end{Theorem}
We will use $K_{\bT}(X)$ to denote the $\bT$-equivariant $K$-theory ring of $X$. Due to the last theorem we do not  distinguish between the algebraic and the topological versions. This theorem implies that:
$$
K_{\bT}(X)=\matZ[ x_{i,j}^{\pm 1}, \aa^{\pm}, \hbar^{\pm 1} ]/R 
$$
where $x_{i,j}$, $i\in I$, $j=1,\dots,\textsf{v}_{i}$ denote the Grothendieck roots of the tautological bundles, $\aa$, $\hbar$ stand for the equivariant parameters of $\bT$ and $R$ denotes certain ideal.  If $X^{\bT}$ is finite, then the ideal $R$ can be described as the ideal of Laurent polynomials whose restrictions (\ref{cgres}) to all fixed points in $X^{\bT}$ vanish.

 As a corollary of Theorem \ref{tmeven}, 
$$\textrm{Pic}(X)\cong \matZ^{|I|}$$
is a finite dimensional lattice generated by the classes of the tautological line bundles 
\be \label{linbundl}
\cL_i=\det(\tb_i)=x_{i,1}\cdots x_{i,\textsf{v}_{i}}, \ \    i \in I.
\ee

\subsection{}
We denote $\bK=\Pic(X)\otimes \matC^{\times}$. The coordinates on torus $\bK$ are called {\it K\"ahler parameters of $X$} (or Fayet–Iliopoulos parameters in physical literature). We use symbols $\zz=(z_i)$, $i\in I$ to denote these characters. 

The coordinates on $\bT=\matC^{\times}_{\hbar}\times \bA$ are called {\it equivariant parameters} (or mass parameters in physical literature) and are denoted by symbols $\hbar$ and $\aa=(a_i)$, $i=1,\dots,\dim(\bA)$.

\subsection{}
Assume $X^{\bT}$ is finite.  The real Lie algebra $\Lie_{\matR}(\bA)$ is naturally equipped with a set of hyperplanes $\{ \omega^{\perp} \}$, where $\omega$ runs over the set of $\bA$-characters appearing as $\bA$-weights of the tangent spaces $T_{p} X$ for $p\in X^{\bT}$. 
The complement of these hyperplanes is a union of connected components 
\be \label{chamb}
\Lie_{\matR}(\bA)\setminus \{ \omega^{\perp}\}= \bigcup \fC
\ee
which are called {\it chambers}.  For a chamber $\fC$ we have the corresponding decomposition
$$
T_p X= N^{+}_p\oplus N^{-}_{p}
$$
where $N^{+}_p,N^{-}_{p}$ are the subspaces whose $\bA$-characters take positive or negative values on $\fC$, respectively. The subspaces $N_p^{\pm}$ are sometimes referred to as {\it attracting} or {\it repelling} directions of the tangent space corresponding to $\fC$.


\subsection{}
Recall that the vertex function ${\bf V}(\aa,\zz)$ of $X$ is a $K$-theoretic analog of Givental's $J$-function in Gromov-Witten theory.
It is defined by the equivariant count of quasimaps from
 a rational curve $\mathbb{P}^1$ to $X$. For the definitions of quasimap moduli spaces and vertex functions see Section 7.2 in \cite{pcmilect}, Section 3 in \cite{AOF} or Section 6 in \cite{AOElliptic}.

 By definition, the vertex function is a power series in the K\"ahler parameters:
 $${\bf V}(\aa,\zz) \in K_{\bT\times \matC^{\times}_q }(X)_{loc}[[\zz]]$$ where $K_{\bT\times \matC^{\times}_q }(X)_{loc}$ denotes the localized equivariant $K$-theory ring of $X$ and $\matC^{\times}_q$ is a one-dimensional torus acting on the domain of the quasimaps $\mathbb{P}^1$ by scaling the homogeneous coordinates  such that the character of $T_{0}\mathbb{P}^1$ is equal to $q$.

The classes of the fixed points $p\in X^{\bT}$ provide a natural basis in localized $K$-theory. The components of the vertex function in this basis have the following form:
\be \label{verp}
{\bf V}_{p}(\aa,\zz) = \sum\limits_{d\in C_{\textrm{eff}}(X)\cap H^{2}(X,\matZ)} c_{d}(\aa) \zz^d
\ee 
where the degrees $d$ run over a certain cone $C_{\textrm{eff}}(X)\subset H^{2}(X,\matR)=\Lie_{\matR}(\bK)$ spanned by effective curves\footnote{The choice of the cone $C_{\textrm{eff}}(X)$ corresponds to the choice of the stability parameter for the Nakajima variety $X$.}, and $c_{d}(\aa) \in \matZ(\aa,q,\hbar)$ are some rational functions.

Alternatively, the vertex function can be defined as the {\it index} of the  $3d$ gauge theory on $\matC\times S^{1}$. 
The choice of the fixed point $p$ in (\ref{verp}) corresponds then to the choice of the vacuum of the theory at infinity of  $\matC\times S^{1}$,
see Section 8.5 in \cite{AOF} for a physical discussion in this direction.

Finally, we also mention that the vertex functions play an important role in the theory of integrable systems known as $XXZ$-spin chains, which were recently understood as the quantum
$K$-theory  of Nakajima varieties, see \cite{NS1,NS2,OkBethe,Pushk1,Pushk2,GaKor}.

\subsection{}
We will be interested in specializations of the vertex functions
corresponding to vanishing equivariant parameters. More precisely, let $\sigma: w\in \matC^{\times} \to \bA $ be a cocharacter from a chamber $\fC$. We define
$$
{\bf V}_{p}(0_{\fC},\zz):=\lim\limits_{w\to 0} {\bf V}_{p}(\sigma(w),z) \in \matZ(q,\hbar)[[\zz]].
$$
It follows from the construction of the virtual structure sheaf on the moduli space of quasimaps that these limits exist and are well defined for all chambers.

\subsection{}
$3d$-mirror symmetry is a conjecture which claims that for every symplectic variety from the class discussed above there exists a dual variety $X'$ called the {\it $3d$-mirror of $X$}. The geometries of $X$ and $X'$ intimately related. 

A $3d$ mirror $X'$  can be characterized in the language of elliptic stable envelopes \cite{MirSym1,MirSym2}\footnote{See also the talk by A.Okounkov ``\textit{Enumerative symplectic duality}'' at the MSRI workshop   \textbf{Structures in Enumerative Geometry} in April 2018 for the first discussion of these ideas (available online). }. The definition in Section 1.1 in \cite{MirSym2} is especially convenient for us.  This duality, among other things, provides the following data:

\begin{itemize}
	
\item An isomorphisms of tori
\be \label{3dmiriso}
\kappa: \, \bA'\times \bK'\times \matC_{\hbar'}^{\times} \times \matC^{\times}_{q} \ \ \longrightarrow \bA \times \bK \times \matC_{\hbar}^{\times} \times \matC^{\times}_{q}
\ee
where $\bK'$ and $\bA'\times \matC_{\hbar'}^{\times}=\bT'$ are the K\"ahler and the equivariant tori of $X'$. In coordinates, the map $\kappa$ is an invertible monomial transformation of the form
\be \label{ident}
z_i \to  {\hbar'}^{n_i} \prod  \limits_{j} {a'}_{j}^{m_{i,j}} ,  \ \ \ a_i \to  ({\hbar'/q})^{r_i} \prod  \limits_{j} {z'}_{j}^{s_{i,j}}, \ \ \hbar\to q/\hbar' 
\ee 
for some $n_i,r_i,m_{i,j},s_{i,j}\in \matZ$.	

\item A bijection of sets of fixed points
\be \label{bij}
\textsf{b}: X^{\bT} \to  (X')^{\bT'}
\ee	 
\end{itemize}
The restriction of (\ref{ident}) to $ker(q)\cap ker(\hbar)$ provides an isomorphisms of tori:
$$
\bar{\kappa}: \bA  \rightarrow \bK', \ \ \ \bK  \rightarrow \bA'
$$
One often refers to this by saying that ``3d mirror symmetry identifies equivariant parameters with K\"ahler parameters of the dual theory''.  
In addition, it provides an identification of chambers with effective cones of the $3d$ mirror variety:
$$
d \bar{\kappa} (\fC)=C_{\textrm{eff}}(X'), \  \ \ \ d \bar{\kappa} (C_{\textrm{eff}}(X))=\fC',
$$
where $d\bar{\kappa}$ is the induced map of the Lie algebras. 
It is therefore more natural to think that the $3d$ mirror symmetry provides a pair $(X',\fC')$ for each pair $(X,\fC)$.

\subsection{}
The isomorphism (\ref{3dmiriso}) provides the dual variety $X'$ with an action of the torus~$\bK$.  Given a fixed point $p\in X^{\bT}$ it is therefore natural to ask 

\begin{center}
{\it What is the $\bK$-character of the tangent space $T_{\textsf{b}(p)} X'$ in terms of $X$? }
\end{center}

The $3d$-mirror of a symplectic variety is not known in most cases. Therefore,
the above question is non-trivial. The answer to this question provides new interesting information about the geometry of $X'$.

In this paper, we propose a combinatorial formula for the character of $T_{\textsf{b}(p)} X'$. The feature of our formula  is that it describes the corresponding character purely in terms of the vertex functions of $X$. In particular, to compute the characters of $T_{\textsf{b}(p)} X'$, one does not need to know what $X'$ is.

\subsection{}
For $N \in K_{\bT'}(pt)$ given by a polynomial $N=w_1+\dots +w_m$ we abbreviate 
$$
\Xi(b,N)=\xi(b,w_1)\cdots \xi(b,w_m),
$$  
where 
\be \label{qbinom}
\xi(b,w)=\dfrac{\varphi(b w)}{\varphi(w)} \ \ \ \textrm{and}  \ \ \  \varphi(w)=\prod\limits_{n=0}^{\infty} (1- w q^n)
\ee
is the $q$-analog of the Gamma function (we assume $0<|q|<1$ so that the last infinite product is convergent).
By the $q$-binomial theorem 
$$
\xi(b,w)=\sum\limits_{n=0}^{\infty}\, \dfrac{(b)_n}{(q)_n} w^n
$$
where
$$
(a)_n = \dfrac{\varphi(a)}{\varphi(a q^n)} = 
\left\{\begin{array}{ll}(1-a) \cdots (1- a q^{n-1}), & n\geq0 \\
\dfrac{1}{(1-a q^{n})\cdots(1-a q^{-1}) }, & n<0\\
\end{array}\right. 
$$
and thus we may think of $\Xi(b,N)$ as a power series in  $w_i$. Similarly, we denote
\be \label{phidef}
\Phi(N)=\varphi(w_1)\dots \varphi(w_m),
\ee
and extend it by multiplicativity to polynomials with negative coefficients. For example:
$$
\Phi(a+2b-3c)=\dfrac{\varphi(a) \varphi(b)^2}{\varphi(c)^3}.
$$

\subsection{} 

Let $(X,\fC)$ and $(X',\fC')$ be pairs of symplectic varieties and chambers related by $3d$-mirror symmetry. 
For $p\in X^{\bT}$ let $\textsf{b}(p)$ be the corresponding fixed point on the dual variety $X'$. The chamber $\fC'=d\bar{\kappa}(C_{\textrm{eff}}(X))$ provides a decomposition
$$
T_{\textsf{b}(p)} X' = N^{+}_{\textsf{b}(p)}\oplus N^{-}_{\textsf{b}(p)}.
$$
We identify these spaces with their $K$-theory classes $N^{
\pm}_{\textsf{b}(p)} \in K_{\bT'}(pt)$ (i.e., simply the $\bT'$-character of this vector space)

\begin{Conjecture}
{\it The vertex functions of $X$ with vanishing equivariant parameters are given by the Taylor series expansions of the following functions
\be\label{hypot}
\begin{array}{|c|}  
\hline \\
\ \ \ \kappa^{*} {\bf V}_p(0_{\fC},\zz)= \Xi(q/\hbar', (N^{-}_{\textsf{b}(p)})^{*} )\ \ \ \\
\\
\hline
\end{array}
\ee 	
where $\kappa^{*}$ stands for substitution (\ref{ident}) and $(N^{-}_{\textsf{b}(p)})^{*}$ is the $\bT'$-module dual to $N^{-}_{\textsf{b}(p)}$ (i.e. the weights of $(N^{-}_{\textsf{b}(p)})^{*}$ are inverses of the weights of $N^{-}_{\textsf{b}(p)}$). }
\end{Conjecture}

Below we check this conjecture in several cases by explicit computation.

\begin{Note}
Note that the right side of (\ref{hypot}) is a product of $\dim(X'/2)$  $q$-binomials (\ref{qbinom}) with $b=\hbar=q/\hbar'$. Thus, the conjecture provides a non-trivial summation formula for the power series defining the vertex functions. 	
\end{Note}

\begin{Note}
Note that the left side of the proposed identity is defined purely in terms of~$X$. Computing the vertex function of $X$
we can therefore obtain the character of $N^{-}_{\textsf{b}(p)}$ and thus the character of the whole tangent space on the dual side:
$$
T_{\textsf{b}(p)} X'=N^{-}_{\textsf{b}(p)} + (N^{-}_{\textsf{b}(p)})^{*}/{\hbar'}.
$$ 
\end{Note}

\section*{Acknowledgments} 
This work is supported by the Russian Science Foundation under grant 19-11-00062.

\section{Vertex functions for Nakajima varieties}
In this section we recall how to compute vertex function (\ref{verp}) for a Nakajima variety $X$. For more details we refer to  Section 7.2 in \cite{pcmilect}, Section 3 in \cite{AOF} or Section 6 in \cite{AOElliptic}.

\subsection{} 
Let $X$ be a Nakajima variety.  We denote by $\textsf{Q}$ be the corresponding quiver and by $I$ the set of vertices of $\textsf{Q}$.  Let  $V_i$ and $W_i$, $i\in I$  be the vector spaces and framing vector spaces corresponding to $X$. We denote $\textsf{v}_{i}=\dim V_i$,  $\textsf{w}_{i}=\dim W_i$.

Let us fix an orientation of $\textsf{Q}$ and define the polarization by the following virtual vector space
$$
P=\Big(\bigoplus_{i\to j} V_i^{*} \otimes V_j\Big) \bigoplus \Big(\bigoplus_{i\in I} W_i^{*} \otimes V_i\Big) - \bigoplus_{i\in I} V_i^{*} \otimes V_i 
$$
where $i\to j$ denotes the sum running over arrows of $\textsf{Q}$ (with fixed orientation). We denote by $\mathscr{P} \in K_{\bT}(X)$ the virtual bundle associated 
to $P$. Explicitly, this K-theory class is given by the following Laurent polynomial:
\begin{align}  \nonumber
\mathscr{P}=\sum\limits_{i\to j} \left(\sum_{m=1}^{\textsf{v}_{i}} x_{i,m}^{-1}\right) \left(\sum_{m=1}^{\textsf{v}_{j}} x_{j,m} \right)+\sum\limits_{i\in I} \left(\sum_{m=1}^{\textsf{w}_{i}} a_{i,m}^{-1}\right) \left(\sum_{m=1}^{\textsf{v}_{i}} x_{i,m} \right)\\ \label{polar}
-\sum\limits_{i\in I} \left(\sum_{m=1}^{\textsf{v}_{i}} x_{i,m}^{-1}\right) \left(\sum_{m=1}^{\textsf{v}_{i}} x_{i,m} \right)
\end{align}
where $x_{i,j}$ and $a_{i,j}$ denote the Grothendieck roots of the tautological bundles $\tb_i$ and the equivariant parameters associated to the framings $W_i$.

\subsection{} 
For the tautological line bundles (\ref{linbundl}) we define the following formal expression in Grothendieck roots
$$
{\bf e}({\bs x},\zz):=\exp\Big(\dfrac{1}{\ln{q}} \sum\limits_{ i \in I } \ln(z_i)  \ln(\cL_{i}) \Big)=\exp\Big(\dfrac{1}{\ln{q}} \sum\limits_{ i \in I }\sum_{j=1}^{\textsf{v}_i} \ln(z_i)  \ln(x_{i,j}) \Big)
$$
where ${\bs x}$ abbreviates the set of all Grothendieck roots. This function has the following transformation properties: 

$$
{\bf e}(x_{1,1},\dots,x_{i,j}q,\dots  ,\zz)=z_i\,{\bf e}(x_{1,1},\dots,x_{i,j},\dots  ,\zz)
$$ 
and similarly
$$
{\bf e}({\bs x},\dots,z_1,\dots, z_{i} q ,\dots, z_{|I|})=\cL_i\,{\bf e}({\bs x},\dots,z_1,\dots, z_{i} ,\dots, z_{|I|}).
$$

\subsection{}
Let $p\in X^{\bT}$ be a fixed point on the quiver variety. The restriction of a tautological bundle $\tb_i$ to a fixed point provides a polynomial:
$$
\left.\tb_i\right|_{p}=t_{i,1}+\dots+t_{i,\dim V_{i}} \in K_{\bT}(pt). 
$$
This polynomial can be defined as ``evaluation'' of the Grothendieck roots  of $\tb_i$ at the corresponding fixed point:
\be \label{cgres}
x_{i,j} \to x_{i,j}(p):=t_{i,j} \in K_{\bT}(pt).
\ee
Note that all $K$-theory classes are symmetric polynomials in $x_{i,j}$ 
and thus the last operation is well defined. We denote by 
${\bs x}(p)$ the set of $x_{i,j}(p)$.

\subsection{} 

Up to normalizing prefactor, the vertex function for Nakajima  varieties equals
\be \label{vertild}
{\bf \tilde V}_p(\aa,\zz)= \prod\limits_{i,j} \int\limits_{0}^{x_{i,j}(p)}  d_q x_{i,j} \, \Phi\Big((q-\hbar) \mathscr{P} \Big) {\bf e}({\bs x},\zz) 
\ee
where the integral symbols denote the Jackson $q$-integral over all Grothendieck roots:
$$
\int\limits_{0}^{a}  d_q x f(x) :=\sum\limits_{n=0}^{\infty}\, f(a q^n),
$$
$\mathscr{P}$ is the polarization (\ref{polar}) expressed in Grothendieck roots and $\Phi$ is defined by (\ref{phidef}). The normalized vertex function equals:
\be \label{verpow}
{\bf V}_p(\aa,\zz)= \frac{1}{\Phi\Big((q-\hbar) \mathscr{P}(p) \Big) {\bf e}({\bs x}(p),\zz)} {\bf \tilde V}_p(\aa,\zz),
\ee
where $\mathscr{P}(p)$ is the result substituting (\ref{cgres}) into the polarization. Thanks to the identity $\varphi(a)/\varphi(a q^d)~=~(a)_{d}$ the function ${\bf V}_p(\aa,\zz)$ is a power series in $\zz$ with coefficients in $\matZ(\aa,\hbar,q)$. The normalization (\ref{verpow}) has the effect of starting the power series with 1:
$$
{\bf V}_p(a,z)= 1+O(\zz) \in \matZ(\aa,\hbar,q)[[\zz]]. 
$$ 	
\begin{Theorem}[Theorem 4 in \cite{AOElliptic}]
{\it For all $p\in X^{\bT}$ the functions (\ref{vertild}) are solutions of certain $q$-difference equations in the equivariant and K\"ahler parameters with regular singularities.}   	
\end{Theorem}
For Nakajima quiver varieties the corresponding $q$-difference equations were described in \cite{OS}.  As a corollary, the power series (\ref{verpow})
converges to an analytical function in some neighborhood $|z_i|<\epsilon_i$ with $0<\epsilon_i$.

In the following sections we compute the vertex functions explicitly for several Nakajima varieties and check the identity (\ref{hypot}).

\section{Cotangent Bundle of $Gr(k,n)$}

\subsection{}
Let $X=T^*Gr(k,n)$, the cotangent bundle of the Grassmannian of $k$-planes in $n$-dimensional space. $X$ is an example of a Nakajima quiver variety, see Section 3 in \cite{MirSym1}. The corresponding quiver $Q$ consists of a single vertex. The associated representation of the quiver consists of a vector space $V$ with $\dim V =k$ and a framing space $W$ with $\dim W = n$. 
Let $\bA=(\matC^{\times})^{n}$ be the torus acting on $\matC^n$ by 
$$
(z_1,\ldots, z_n)\to (z_1 a_1^{-1}, \ldots, z_n a_n^{-1}).
$$
$\bA$ induces a natural action on $X$. Let $\matC^{\times}_{h}$ be the torus which acts on $X$ by 
$$
\eta_p \to \hbar^{-1} \eta_p
$$
where $\eta_p$ is a covector at a point $p\in X$. Thus, the torus $\bT=\bA \times \matC_{\hbar}^{\times}$ acts on~$X$. One can check that $\bA$ preserves the sympectic form, while $\matC_{\hbar}^{\times}$ scales it by~$\hbar$.

\subsection{}

The set of fixed points $X^{\bT}$ corresponds to the $k$-dimensional coordinate subspaces in $X$, and thus fixed points $\textbf{p}\in X^{\bT}$ are naturally labeled by size $k$ subsets of $\{1, \ldots, n\}$. As a Nakajima quiver variety, $X$ comes equipped with a rank $k$-tautological vector bundle $\tb$. Let $\textbf{p}=\{ p_1 , \ldots , p_k \}\subset \{1, \ldots, n\}$. The $\bT$-character of the $p_i$-th coordinate direction is $a_{p_i}^{-1}$. Thus the restriction of the tautological bundle to the fixed point $\textbf{p}$ is given by 
$$ 
\tb_{\textbf{p}} = \sum_{i=1}^k a_{p_i}^{-1}.
$$
Thus, the restriction of the Grothendieck roots in (\ref{cgres}) has the form:
\be \label{fxdpt}
x_i(\textbf{p}) = a_{p_i}^{-1}.
\ee

\subsection{}
In this case, the polarization (\ref{polar}) takes the form 
$$
P=  W^* \otimes V - V^* \otimes V 
$$
and we have 
$$
\Phi((q-h)\cP)= \prod_{i=1}^n \prod_{j=1}^k \frac{\phi(q a_i x_j)}{\phi(h a_i x_j)} \prod_{i,j=1}^k \frac{\phi(h x_i/x_j)}{\phi(qx_i/x_j)}
$$
Thus, the vertex functions associated to $T^* Gr(k,n)$ are given by the power series:
$$
\begin{aligned}
& {\bf V}_{\textbf{p}}(a_1, \ldots, a_n, z) = & \\ & \sum_{d_{1}, \ldots, d_{k}=0}^{\infty} \left(\prod_{i=1}^n \prod_{j=1}^k \frac{(h a_i /a_{p_j} )_{d_{j}}}{(q a_i /a_{p_j})_{d_{j}}} \prod_{i,j=1}^k \frac{ (q a_{p_i}/a_{p_j} )_{d_i-d_j}}{(h a_{p_i} /a_{p_j})_{d_{i}-d_{j}}} \right) z^{d_{1}+\ldots + d_{k}}. 
\end{aligned}
$$

\begin{Note}
In the case of $k=1$, the fixed point $\bf p$ is just an integer $1 \leq p \leq n$, and the vertex function is 
$$
  {\bf V}_{\textbf{p}}(a_1, \ldots, a_n, z) =  \sum_{d=0}^{\infty} \left(\prod_{i=1}^n  \frac{(h a_i /a_{p} )_{d}}{(q a_i /a_{p})_{d}}  \right) z^{d}.
$$
which is an example of the well-known $q$-hypergeometric series (see \cite{hypergeo}).
\end{Note}

\subsection{}

The torus $\bA$ is $n$-dimensional, and we choose the chamber $\fC$ associated to the cocharacter $w \mapsto (w^{-1}, w^{-2}, \ldots, w^{-n})$. Then
$$
{\bf V}_{\bf p}(0_{\fC}, z)= \lim_{w \to 0} {\bf V}_{\bf p} (w^{-1}, \ldots, w^{-n}, z)
$$
In terms of equivariant parameters, this limit corresponds to $\frac{a_i}{a_j} \to 0$ if $i<j$.
Substituting (\ref{fxdpt}) and taking the limit, we compute 
$$
{\bf V}_{\bf p}(0_{\fC}, z) = \sum_{d_{1}, \ldots , d_{k}=0}^{\infty} \left( \prod_{i=1}^k \frac{(\hbar)_{d_i}}{(q)_{d_i}} \prod_{j=i+1}^{n} \left( \frac{\hbar}{q} \right)^{d_i} \prod_{i=1 }^k \prod_{j=i+1}^k \left( \frac{q}{\hbar} \right)^{d_i-d_j} \right) z^{d_1+\ldots+ d_k} 
$$
which simplifies to 
$$
{\bf V}_{\bf p}(0_{\fC}, z) = \sum_{d_1, \ldots, d_k=0}^{\infty} \left( \prod_{i =1}^{k} \frac{(\hbar)_{d_i}}{(q)_{d_i}} \left( \left(\frac{\hbar}{q} \right)^{n-p_i-k+2i-1} z \right)^{d_i} \right).
$$
Interchanging the order of the summation and the product and using the $q$-binomial theorem, we obtain
\be \label{limvergr}
{\bf V}_{\bf p}(0_{\fC}, z)=\prod_{i =1}^{k} \xi
\left(\hbar, \left(\frac{\hbar}{q}\right) ^{n-k-p_i+2i-1} z \right) .
\ee

\subsection{}
For $n\geq 2k$ the 3d mirror dual $X'$ of $T^* Gr(k,n)$ is the Nakajima quiver variety associated to the $A_{n-1}$ quiver with dimension vector 
$$
\mathsf{v} = (1, 2 , \ldots, k-1, k, \ldots, k , k-1, \ldots , 2 ,1)
$$
with one dimensional framings at positions $k$ and $n-k$ and trivial framings elsewhere (see Section 4 of \cite{MirSym1} for details). The fixed points of $X'$ are labeled by Young diagrams $\lambda$ which fit into the rectangle with dimensions $(n-k)\times k$. The boundary of such a diagram is a piecewise linear curve with $n-k$ horizontal segments and $k$ vertical segments. The positions of the $k$ vertical segments naturally determine a $k$-tuple ${\bf p}_{\lambda} \subset \{1,2,\ldots,n\}$. The bijection on fixed points (\ref{bij}) has the form (see also Section 6.1 in \cite{MirSym1}):
$$
\textsf{b}: \lambda \to {\bf p}_{\lambda}
$$
The isomorphism (\ref{ident}) between the equivariant parameters of $X'$ and the K\"ahler parameters of $X$ is the following:
$$
\kappa^{*}: z\to a' {\hbar'}^{k-1}, \ \ \hbar \to q/\hbar'
$$
where $a'=a'_1/a'_2$ is the equivariant parameters of $X'$ (i.e. coordinates on the framing torus $\bA' = (\matC^{\times})^2$). 
Using this transformation we can rewrite (\ref{limvergr}) as follows:
$$
\kappa^{*} \Big({\bf V}_{\bf p}(0_{\fC}, z)\Big)=\Xi \Big(q/\hbar',\sum\limits_{i=1}^{k} a' {\hbar^{'}}^{2k-n+p_i-2 i}\Big).
$$

We choose the chamber $\fC'$ given by the cocharacter $w \mapsto (a_1',a_2')= (1,w^{-1})$.
The characters of the tangent spaces at the fixed points and of repelling and attracting parts of it corresponding to $\fC'$ were computed in Section 4.4 of \cite{MirSym1}:
$$
N'^{+}_{\lambda} = \sum_{i=1}^k \frac{a'_1}{a'_2} {\hbar'}^{2k-n+\textbf{p}_i-2i-1}, \ \ N'^{-}_{\lambda}= \sum_{i=1}^k \frac{a'_2}{a'_1} {\hbar'}^{-2k+n-\textbf{p}_i+2i},
$$
where $\textbf{p}=\textsf{b}(\lambda)$. Thus, the character of the dual representation equals:
$$
(N'^{-}_{\lambda})^{*}= \sum_{i=1}^k a' {\hbar'}^{2 k-n+\textbf{p}_i-2i}
$$
We, therefore have:
$$
\kappa^{*} \Big({\bf V}_{\bf p}(0_{\fC}, z)\Big)=\Xi(q/\hbar',(N'^{-}_{\lambda})^{*} ).
$$
which proves Conjecture 1 in this case.


\begin{Note}
As we mentioned above, the dual variety $X'$ is a Nakajima variety if $n\geq 2k$. In the case 	$n<2k$ the variety $X'$ is not known to us. Conjecture~1, however, still holds in this case. 

For example, our formula (\ref{limvergr}) implies that in the case of the zero-dimensional Nakajima variety  $X=T^{*}Gr(n,n)$, the dual variety $X'\cong\matC^{2n}$ is a $2n$-dimensional symplectic vector space  equipped with an action of two-dimensional torus $\bT'=\matC_{a'}^{\times}\times \matC_{\hbar'}^{\times}$, so that
$$
\textrm{char}_{\bT'}(\matC^{2n})=\sum\limits_{i=1}^{n}\,( a'{\hbar'}^{n-i-1} +  {\hbar'}^{-n+i}/a') \in K_{\bT'}(pt).
$$
We believe that the $3d$ mirrors can be constructed as \textit{bow varieties} \cite{Cherk,NakBow}. 
\end{Note}

\section{Hilbert scheme of points on $\matC^2$ \label{hilbsec}} 

\subsection{\label{torin}} 
In this section, $\hilb$ denotes the Hilbert scheme of $n$ points on 
the complex plane. This is a smooth, symplectic, quasiprojective variety which parametrizes 
the polynomial ideals of codimension $n$ \cite{NakajimaLectures1,NakajimaLectures2}:   
$$
\hilb=\{ {\cal{J}} \subset \matC[x,y] :\,  \dim_{\matC} (\matC[x,y]/ {\cal{J}}) =n \}. 
$$	
Let $\bT\cong (\matC^{\times})^2$ be a two-dimensional torus acting on $\matC[x,y]$ by:
\be
\label{polact}
(x,y)\to (x t_1^{-1},y t_2^{-1}).
\ee
This action induces an action of $\bT$ on  $\hilb$. The one-dimensional space spanned by a symplectic form $ \omega \subset H^{2}(\hilb,\matC)$ is a natural $\bT$-module. We denote by $\hbar$ the  $\bT$-character of $\matC \omega$. From our normalization (\ref{polact}) we find:
$$
\hbar^{-1} = t_1 t_2 \in K_{\bT}(pt)=\matZ[t_1^{\pm 1}, t_2^{\pm 1}].
$$
We denote by 
$$
\bA=\ker(\hbar) \subset \bT 
$$
the one-dimensional subtorus preserving the symplectic form on $\hilb$. We denote the coordinate on $\bA$ by $a$, such that\footnote{We may assume $\hbar^{1/2}$ exists by passing to the double cover of $\bT$ if needed. }:
\be \label{ttoa}
t_1 =\dfrac{a}{\sqrt{\hbar}}, \ \ t_2=\dfrac{1}{a \sqrt{\hbar}}.
\ee

\subsection{}
The set of fixed points $\hilb^{\bT}=\hilb^{\bA}$ is a finite set labeled by Young diagrams with $n$ boxes. 
A Young diagram $\lambda$ corresponds to the $\bT$-invariant ideal generated by monomials:
$${\cal{J}}_{\lambda}=\{x^{\lambda_1}, x^{\lambda_1} y, x^{\lambda_2} y^2 \dots\}.$$
The Hilbert scheme is equipped with a rank $n$-tautological vector bundle $\tb$ with fibers at a point ${\cal{J}}\in X$ given by
$$
\left.\tb\right|_{\cal{J}}=\matC[x,y]/ {\cal{J}}.
$$
Let $\Box=(m,n)$ be a box with coordinates $m$ and $n$ in a Young diagram~$\lambda$. The $\bT$-character of the corresponding monomial equals
$
t_1^{-(n-1)} t_2^{-(m-1)}.
$
We conclude that restriction of the tautological bundle to the fixed point $\lambda$ is given by the sum of the characters of monomials over the boxes in $\lambda$:
$$
\left.\tb\right|_{\cal{J}_{\lambda}}=\sum\limits_{(m,n)\in \lambda} t_1^{-(n-1)} t_2^{-(m-1)}.
$$
Thus, the restriction of Grothendieck roots in (\ref{cgres}) has the form:
$$
x_i(\lambda) =  t_1^{-(n_i-1)} t_2^{-(m_i-1)},
$$
where $(m_i,n_i)$ are the coordinates of the $i$-th box in $\lambda$. 
Again, the order of boxes is not important, because all tautological classes are symmetric in $x_i$'s. 

\subsection{} 
The characters of the tangent space at a fixed point $\lambda$ can be computed using, for example, a free resolution of the ideal ${\cal{J}}_{\lambda}$, see Section 3.4.19 \cite{pcmilect}. This gives the following explicit formula:
$$
T_{\lambda} X = \sum\limits_{i\in \lambda} \, t_1^{-l_{\lambda}(i)} t_2^{a_{\lambda}(i)+1}+t_2^{-a_{\lambda}(i)} t_1^{l_{\lambda}(i)+1} \in K_{\bT}(pt),
$$
where the sum runs over the boxes of the Young diagram $\lambda$ and 
$$
l_{\lambda}(\Box)=\lambda_i-j, \ \  a_{\lambda}(\Box)=\lambda^{'}_j-i
$$
stand for the standard leg and arm length of a box $i$ in $\lambda$ ($\lambda'$ denotes the transposition of $\lambda$).  The change of variables (\ref{ttoa}) provides a decomposition of this character into a sum of characters with positive and negative powers of $a$:
$$
T_{\lambda} X = N_{\lambda}^{+} + N_{\lambda}^{-}
$$ 
where 
\be \label{negpart}
N_{\lambda}^{+}=\sum\limits_{i\in \lambda} \, t_2^{-a_{\lambda}(i)} t_1^{l_{\lambda}(i)+1},  \ \ \ N_{\lambda}^{-}=\sum\limits_{i\in \lambda} \, t_1^{-l_{\lambda}(i)} t_2^{a_{\lambda}(i)+1}.
\ee

\subsection{} 
The Hilbert scheme is an example of a Nakajima quiver variety, see Section 3.3 in  \cite{SmirnovElliptic}. 
The corresponding quiver $\textsf{Q}$ consists of one vertex and one loop
connecting the vertex with itself. The associated representation of the quiver consists of a vector space $V$ with
$\dim V=n$ and a framing space $W$ with $\dim W=1$. 

The polarization (\ref{polar}) in this case has the form (see Section 3.4 in \cite{SmirnovElliptic}):
$$
P=W^{*} \otimes V + V^{*} \otimes V t_1 - V^{*} \otimes V
$$
We compute:
$$
\Phi((q-\hbar) \cP)=\prod\limits_{i=1}^{n} \dfrac{\varphi(q x_i)}{\varphi(\hbar x_i)} \prod\limits_{i,j=1}^{n} \dfrac{\varphi(q t_1 x_i/x_j )}{\varphi(\hbar t_1 x_i/x_j )} \dfrac{\varphi(\hbar x_i/x_j )}{\varphi(q x_i/x_j )}. 
$$
Thus, the vertex functions associated to the Hilbert scheme $X$ are given by the following power series: 
$$
\begin{aligned} \label{verhilb}
& {\bf V}_{\lambda}(a,z)=  \\
& \sum\limits_{d_1,\dots,d_n=0}^{\infty} \, \Big(\prod\limits_{i=1}^{n} \dfrac{(\hbar x_i(\lambda))_{d_{i}}}{(q x_i(\lambda))_{d_{i}}} \prod\limits_{i,j=1}^{n} \dfrac{(\hbar t_1 x_i(\lambda)/x_j(\lambda) )_{d_i-d_j}}{(q t_1 x_i(\lambda)/x_j(\lambda) )_{d_i-d_j}} \dfrac{(q x_i(\lambda)/x_j(\lambda) )_{d_i-d_j}}{(\hbar x_i(\lambda)/x_j(\lambda) )_{d_i-d_j}} \Big) z^{d_1+\dots +d_n} 
\end{aligned}
$$
which is a power series in $z$ with coefficients in $\matZ(\hbar,q,a)$. 

\subsection{}
In the case of the Hilbert scheme $X$, the torus $\bA$ preserving the symplectic form is one-dimensional and the chamber decomposition (\ref{chamb}) takes the form:
$$
\Lie_{\matR}(\bA)\setminus \{\omega^{\perp}\}= \matR\setminus \{0\} =\fC_{+}\cup \fC_{-}.  
$$ 
Let us analyze the case $\fC=\fC_{+}$ which corresponds to the cocharacters 
$\{a\to 0\}$ so that 
\be \label{limver}
{\bf V}_{\lambda}(0_{\fC},z)=\lim\limits_{a\to 0} {\bf V}_{\lambda}(a,z).
\ee
It is expected that the Hilbert scheme $X$  is self-dual with respect to $3d$ mirror symmetry\footnote{To the best of authors knowledge this is not proved.}:
$$
X^{'}\cong X.
$$
We denote by $t_1',t_2',\hbar',a'$ the equivariant parameters of $X'$. As in (\ref{ttoa}) these parameters are related by
$$
t_1'=\dfrac{a'}{\sqrt{\hbar'}}, \ \ \ t_1'=\dfrac{1}{a' \sqrt{\hbar'}}.
$$
We consider the trivial bijection on the fixed points (\ref{bij})
$$
\mathsf{b} = id: \lambda \mapsto \lambda,
$$
and the change of variables (\ref{ident}):
\be \label{varch}
\kappa^{*}: z\to a' \sqrt{\hbar'}, \ \ \hbar \to \dfrac{q}{\hbar'}.
\ee
As above, let $\kappa^{*} {\bf V}_{\lambda}(0_{\fC},z)$ denote the result of the substitution (\ref{varch}) into the power series given by the limit (\ref{limver}). Thus, $\kappa^{*} {\bf V}_{\lambda}(0_{\fC},z)$ is a power series in~$a'$:
$$
\kappa^{*} {\bf V}(0_{\fC},z)\in \matZ(\sqrt{\hbar'},q)[[a']].  
$$ 
 Then, the equality (\ref{hypot}) in this case is equivalent to the following identity
\be
\begin{array}{|c|}
	\hline\\
\kappa^{*} {\bf V}(0_{\fC},z)=\Xi(q/\hbar^{'},(N^{-}_{\lambda})^{*}) = \prod\limits_{i\in \lambda}\, \xi(q/\hbar^{'}, {t_1'}^{l_{\lambda}(i)} {t_2'}^{-a_{\lambda}(i)-1})\\
\\
\hline
\end{array} 
\ee
where $N^{-}_{\lambda}$ is from (\ref{negpart}). 
This is a quite nontrivial summation formula for the power series ${\bf V}_{\lambda}(0_{\fC},z)$.
This identity can be proved by induction on the number of boxes. We plan to give a general proof of such identities arising for $A_n$ and affine $\widehat{A}_n$ quiver varieties in the sequel paper.

\bibliographystyle{alpha}
\bibliography{bib}

\def\cprime{$'$} \def\cprime{$'$}
\begin{thebibliography}{RSVZ19b}

\bibitem[AFO18]{AOF}
Mina Aganagic, Edward Frenkel, and Andrei Okounkov.
\newblock {Quantum $q$-Langlands Correspondence}.
\newblock {\em Trans. Moscow Math. Soc.}, 79:1--83, 2018.

\bibitem[AO16]{AOElliptic}
Mina {Aganagic} and Andrei {Okounkov}.
\newblock {Elliptic stable envelopes}.
\newblock {\em arXiv e-prints}, page arXiv:1604.00423, Apr 2016.

\bibitem[AO17]{OkBethe}
Mina Aganagic and Andrei Okounkov.
\newblock {Quasimap counts and Bethe eigenfunctions}.
\newblock {\em Moscow Math. J.}, 17(4):565--600, 2017.

\bibitem[BDG15]{BDG}
Mathew {Bullimore}, Tudor {Dimofte}, and Davide {Gaiotto}.
\newblock {The Coulomb Branch of 3d $\mathcal{N}=4$ Theories}.
\newblock {\em arXiv e-prints}, page arXiv:1503.04817, Mar 2015.

\bibitem[{Che}11]{Cherk}
Sergey~A. {Cherkis}.
\newblock {Instantons on Gravitons}.
\newblock {\em Communications in Mathematical Physics}, 306(2):449--483, Sep
  2011.

\bibitem[Gin12]{GinzburgLectures}
Victor Ginzburg.
\newblock Lectures on {N}akajima's quiver varieties.
\newblock In {\em Geometric methods in representation theory. {I}}, volume~24
  of {\em S\'emin. Congr.}, pages 145--219. Soc. Math. France, Paris, 2012.

\bibitem[GK13]{GaKor}
Davide Gaiotto and Peter Koroteev.
\newblock {On Three Dimensional Quiver Gauge Theories and Integrability}.
\newblock {\em JHEP}, 05:126, 2013.

\bibitem[GR90]{hypergeo}
George {Gasper} and Mizan {Rahman}.
\newblock {\em Basic Hypergeometric Series}.
\newblock Cambridge University Press, 1990.

\bibitem[GW09]{GW}
Davide {Gaiotto} and Edward {Witten}.
\newblock {S-Duality of Boundary Conditions In N=4 Super Yang-Mills Theory}.
\newblock {\em Advances in Theoretical and Mathematical Physics}, 13:721--896,
  Jan 2009.

\bibitem[HW97]{HW}
Amihay {Hanany} and Edward {Witten}.
\newblock {Type IIB superstrings, BPS monopoles, and three-dimensional gauge
  dynamics}.
\newblock {\em Nuclear Physics B}, 492(1):152--190, May 1997.

\bibitem[IS96]{SI}
K.~{Intriligator} and N.~{Seiberg}.
\newblock {Mirror symmetry in three dimensional gauge theories}.
\newblock {\em Physics Letters B}, 387(3):513--519, Feb 1996.

\bibitem[Kal09]{Kaledin}
D.~Kaledin.
\newblock Geometry and topology of symplectic resolutions.
\newblock In {\em Algebraic geometry---{S}eattle 2005. {P}art 2}, volume~80 of
  {\em Proc. Sympos. Pure Math.}, pages 595--628. Amer. Math. Soc., Providence,
  RI, 2009.

\bibitem[KPSZ17]{Pushk2}
Peter {Koroteev}, Petr~P. {Pushkar}, Andrey {Smirnov}, and Anton~M. {Zeitlin}.
\newblock {Quantum K-theory of Quiver Varieties and Many-Body Systems}.
\newblock {\em arXiv e-prints}, page arXiv:1705.10419, May 2017.

\bibitem[MN18]{kirv}
Kevin McGerty and Thomas Nevins.
\newblock Kirwan surjectivity for quiver varieties.
\newblock {\em Invent. Math.}, 212(1):161--187, 2018.

\bibitem[MO12]{MO}
Davesh {Maulik} and Andrei {Okounkov}.
\newblock {Quantum Groups and Quantum Cohomology}.
\newblock {\em arXiv e-prints}, page arXiv:1211.1287, Nov 2012.

\bibitem[Nak98]{Nak1}
Hiraku Nakajima.
\newblock Quiver varieties and {K}ac-{M}oody algebras.
\newblock {\em Duke Math. J.}, 91(3):515--560, 1998.

\bibitem[Nak99]{NakajimaLectures1}
Hiraku Nakajima.
\newblock {\em Lectures on {H}ilbert schemes of points on surfaces}, volume~18
  of {\em University Lecture Series}.
\newblock American Mathematical Society, Providence, RI, 1999.

\bibitem[Nak16]{NakajimaLectures2}
Hiraku Nakajima.
\newblock More lectures on {H}ilbert schemes of points on surfaces.
\newblock In {\em Development of moduli theory---{K}yoto 2013}, volume~69 of
  {\em Adv. Stud. Pure Math.}, pages 173--205. Math. Soc. Japan, [Tokyo], 2016.

\bibitem[NS09a]{NS2}
Nikita~A. Nekrasov and Samson~L. Shatashvili.
\newblock {Quantum integrability and supersymmetric vacua}.
\newblock {\em Prog. Theor. Phys. Suppl.}, 177:105--119, 2009.

\bibitem[NS09b]{NS1}
Nikita~A. Nekrasov and Samson~L. Shatashvili.
\newblock Supersymmetric vacua and {B}ethe ansatz.
\newblock {\em Nuclear Phys. B Proc. Suppl.}, 192/193:91--112, 2009.

\bibitem[NT16]{NakBow}
Hiraku {Nakajima} and Yuuya {Takayama}.
\newblock {Cherkis bow varieties and Coulomb branches of quiver gauge theories
  of affine type $A$}.
\newblock {\em arXiv e-prints}, page arXiv:1606.02002, Jun 2016.

\bibitem[Oko15]{pcmilect}
Andrei Okounkov.
\newblock {Lectures on K-theoretic computations in enumerative geometry}.
\newblock {\em ArXiv: 1512.07363}, 2015.

\bibitem[OS16]{OS}
Andrei Okounkov and Andrey Smirnov.
\newblock {Quantum difference equation for Nakajima varieties}.
\newblock {\em ArXiv: 1602.09007}, 2016.

\bibitem[PSZ16]{Pushk1}
Petr~P. {Pushkar}, Andrey {Smirnov}, and Anton~M. {Zeitlin}.
\newblock {Baxter Q-operator from quantum K-theory}.
\newblock {\em arXiv e-prints}, page arXiv:1612.08723, Dec 2016.

\bibitem[RSVZ19a]{MirSym2}
R.~{Rim{\'a}nyi}, A.~{Smirnov}, A.~{Varchenko}, and Z.~{Zhou}.
\newblock {Three dimensional mirror self-symmetry of the cotangent bundle of
  the full flag variety}.
\newblock {\em arXiv e-prints}, page arXiv:1906.00134, May 2019.

\bibitem[RSVZ19b]{MirSym1}
Rich{\'a}rd {Rim{\'a}nyi}, Andrey {Smirnov}, Alexand {Varchenko}, and Zijun
  {Zhou}.
\newblock {3d Mirror Symmetry and Elliptic Stable Envelopes}.
\newblock {\em arXiv e-prints}, page arXiv:1902.03677, Feb 2019.

\bibitem[{Smi}18]{SmirnovElliptic}
Andrey {Smirnov}.
\newblock {Elliptic stable envelope for Hilbert scheme of points in the plane}.
\newblock {\em arXiv e-prints}, page arXiv:1804.08779, Apr 2018.

\end{thebibliography}

\newpage

\vspace{12 mm}

\noindent
Hunter Dinkins\\
Department of Mathematics,\\
University of North Carolina at Chapel Hill,\\
Chapel Hill, NC 27599-3250, USA\\
hdinkins@live.unc.edu

\vspace{3 mm}

\noindent
Andrey Smirnov\\
Department of Mathematics,\\
University of North Carolina at Chapel Hill,\\
Chapel Hill, NC 27599-3250, USA;\\
Steklov Mathematical Institute\\
of Russian Academy of Sciences,\\
Gubkina str. 8, Moscow, 119991, Russia

\end{document}